\documentclass[10pt]{article}
\usepackage{amssymb,units,amsmath}
\usepackage{bbold} 
\usepackage{bookman}
\usepackage{graphicx}
\usepackage{dcolumn}
\usepackage{bm,bbm,mathtools,esvect}
\usepackage{slashed,fullpage}
\usepackage{comment}
\usepackage{nicefrac}
\title{  Differential geometry using quaternions } 
\author{ Sergio Giardino\footnote{\tt sergio.giardino@ufrgs.br}\\ 
\footnotesize{\it Departamento de Matem\'atica Pura e Aplicada} \\ 
\footnotesize{\it Universidade Federal do Rio Grande do Sul (UFRGS)}\\ 
\footnotesize{\it caixa postal 15080, 91501-970 Porto Alegre RS}\\
\footnotesize{\it Brazil}}

\begin{document}
\date{}

\newtheorem{theorem}{Theorem}[section] 
\newtheorem{remark}{Remark}[section] 
\newtheorem{lemma}{Lemma}[section] 
\newtheorem{proposition}{Proposition}[section] 
\newtheorem{corollary}{Corollary}[section] 
\newtheorem{definition}{Definition}[section]

\maketitle

\begin{abstract} \noindent 
This paper establishes the basis of the quaternionic differential geometry ($\mathbbm H$DG) initiated in a previous article. 
The usual concepts of curves and surfaces are generalized to quaternionic constraints, as well as the curvature and torsion concepts, differential forms, directional derivatives and the structural equations. The analogy between the quaternionic and the real geometries is obtained using a matrix representation of quaternions. The results evidences the  quaternionic formalism as a suitable language to differential geometry that  can be useful in various directions of future investigation.
\vspace{2mm}

\noindent{\bf keywords:} Curves in Euclidean space; Other special differential geometries; Quaternion and other division algebras.

\vspace{1mm}

\noindent{\bf MSC class. codes:} 53A04; 53A40; 11R52.
\end{abstract}

\tableofcontents
\section{INTRODUCTION} 
Following the approach to quaternionic differential geometry ($\mathbbm H$DG) for curves introduced in \cite{Giardino:2021onv}, this article contains a general framework that encompasses curves and surfaces within the concept of regular constraint that will be defined in a while. 
However, the investigation of quaternionic and  hyper-complex geometries is not introduced in the research contained in this article, which is an element of
the  interest in the interplay between differential geometry, Clifford algebras, and their representations. 

The investigation on quaternionic curves seems to have been initiated  by \cite{Bharathi:1987kbn}, theoretically developed by \cite{Gunes:1994mmt,Girard:2015lnc,Aksoyak:2019}, and  applied on various sorts of curves \cite{Coken:2004otq,Hacisalihoglu:2011gok,Gungor:2011scq,Ilarslan:2013bma,EJPAM1979,Bektas:2016qoc,Kisi:2017aap,Coken:2017nqc,Karadag:2019scq,Kizilay:2020eqc,Kahraman:2020dqn,Gungor:2013qie,Onder:2020gqi,Altun:2017scy,Senyurt:2017scy,Hou:2018gia}, as well as surfaces \cite{Yayli:2016sqc,Yayli:2016csq,Yayli:2017qso,Gok:2017sni,Kocakusakli:2017cks,Karakus:2019ame,Yayli:2019zco,Bekar:2021rsc,Tuncer:2022mgt}, and hyper-surfaces in the quaternionic projective space \cite{Berndt:1991jra,Suh:1997dga}. These investigations are not related to complex manifolds, that comprise a vast field of research in differential geometry, and we mention \cite{Gori:2021ana} as a recent development with various classical references. In fact, this research is a development of the real differential geometry whose relation to the territory of complex manifolds are completely unknown, comprising an exciting direction for future research.

The approach to $\mathbbm H$DG presented in this article is somewhat simpler and  more general, capable of describing three kinds of quaternionic constraints: The one-dimensional constraint is equivalent to a curve, the two-dimensional constraint to a surface, and  a three-dimensional quaternionic constraint, without a real counterpart in three dimensions. Following a way that runs parallel to the real differential geometry, as our first task we define quaternionic Frenet-Serret equations. Before doing that, the next section contains the quaternionic notation that will be used throughout the article.

\section{QUATERNIONS}
This section collects essential facts, and establishes the notation adopted in the article, and it is not aimed as an introduction to the subject, that can be found from various sources, {\it e.g.} \cite{Ward:1997qcn,Garling:2011zz,Rocha:2013qtt,Morais:2014rqc}. We also point out that this section does not contain novel result on quaternions, it simply organizes the most important facts to be used in this paper, and establishes the notation.
However, to the best of our knowledge, the discussion about the polar and symplectic notation is original. Although this section is somewhat long, it will be necessary in order to set quaternions within a convenient form to our purposes.

We define quaternions ($\mathbbm H$) as hyper-complex numbers that, $\,\forall\,q\in\mathbbm H$, it holds that
\begin{equation}\label{hn01}
q=x_0+x_1i+x_2 j+x_3k,\qquad\mbox{where}\qquad x_\mu\in\mathbbm R\qquad\mbox{for}\qquad\mu\in\{0,\,1,\,2,\,3\}.
\end{equation} 
The $\,i,\,j\,$ and $\,k\,$  anti-commuting imaginary units  satisfy
\begin{equation}\label{hn02}
i^2\,=\,j^2\,=\,k^2\,=-1\qquad\mbox{and}\qquad ijk=-1.
\end{equation}
Moreover, 
\begin{equation}\label{hn002}
x_0\qquad\mbox{and}\qquad\vec x=x_1i+x_2 j+x_3k
\end{equation}
are respectively named the scalar, or temporal, and the vector, or spatial, components of the quaternion. {\bf If the scalar component of the quaternion is zero, namely $x_0=0$, the quaternion is called a pure imaginary quaternion.} The multiplication of quaternionic imaginary units satisfies the general rule 
\begin{equation}\label{hn003}
 e_m e_n=-\delta_{mn}+\epsilon_{mnl}e_l,\qquad\mbox{where}\qquad m,\,n,\,l\in\{1,\,2,\,3\},
\end{equation}
$\epsilon_{mnl}$ is the anti-commuting Levi-Civit\`a symbol, and we choose
\begin{equation}\label{hn004}
 e_1=i,\qquad e_2=j,\qquad e_3=k.
\end{equation}
Defining $\,e_0=1\,$, we identify
\begin{equation}\label{hn005}
 e_\mu=\big\{e_0,\,e_m\big\}.
\end{equation} 
as the natural basis to the real linear space defined by quaternions.
In the same fashion as complex numbers, the quaternionic conjugate $\,\overline q,\,$ and the quaternionic norm $\,|q|\,$ are as follows
\begin{equation}\label{hn03}
\overline q\,=\,x_0-\vec x,\qquad\qquad |q|^2\,=\,q\overline q\,=\,x_0^2+|\,\vec x\,|^2,\qquad\mbox{and}\qquad |\,\vec x\,|^2=x_1^2+x_2^2+x_3^2.
\end{equation}
According to Adolf Hurwitz \cite{Hurwitz:1898hqt}, quaternions comprise one of the four division algebras, the other being the reals ($\mathbbm R$), the complexes ($\mathbbm C$) and the octonions ($\mathbbm O$).  Notation (\ref{hn01}) can be named  either the extended or cartesian notation, which alternatively may be written as
\begin{equation}\label{hn04}
q=x_0+\omega |\,\vec x\,|,\qquad \mbox{with}    \qquad\omega=\frac{\vec x}{\left|\,\vec x\,\right|}.
\end{equation}
The unitary imaginary function $\omega$, so that $\omega^2=-1$, commutes with the scalar component $\,x_0\,$ of the quaternion, but every quaternionic number has their own imaginary unit $\omega$ that, in general, neither commutes nor anti-commutes with the quaternionic unit of a different quaternion. If necessary, one may write $\,q=x_0+\omega_x |\,\vec x\,|$ to avoid confusion. 
From \cite{Harvey:1990sca}, we adopt the scalar product for quaternions
\begin{eqnarray}\label{hn10}
\nonumber\big\langle p,\,q\big\rangle &=&\mathfrak{Re}\big[\,p\overline q\,\big]\\
&=&\frac{p\overline q+q\overline p}{2},
\end{eqnarray}
which permit us to define that $p$ and $q$ are orthogonal if
\begin{equation}\label{hn11}
\big\langle p,\,q\big\rangle\,=\,0.
\end{equation}
Immediately, 
\begin{equation}\label{hn006}
 \big\langle e_\mu,\,e_\nu\big\rangle=\delta_{\mu\nu},
\end{equation}
and a  four dimensional real vector space appears after using the inner product (\ref{hn10}) together with the natural  basis (\ref{hn005}). An arbitrary quaternion $\,q\,$ defines a new basis to $\mathbbm H$, where  the four original orthogonal directions, {\bf namely the real direction $e_0=1$ and the three imaginary directions $e_m$, generate a novel basis, such as
\begin{equation}\label{hn007}
 q\mapsto\big\{q,\,e_m q\big\}.
\end{equation}
The orthogonality of the four novel directions can be seen from
\begin{equation}
 \big\langle e_\mu q,\,e_\nu q\big\rangle=\delta_{\mu\nu}|q|^2.
\end{equation}
}
Another fundamental property is the quaternionic parallelism: quaternions $p$ and $q$ are parallel if
\begin{equation}\label{hn12}
\big\langle p,\,q\big\rangle\,=\,\,p\overline q,
\end{equation}
Of course, if $\alpha_0\in\mathbbm R$, thus $\,e_\mu\,$ and $\,\alpha_0 e_\mu,\,$  are parallel, as well as $q$ and $\alpha_0 q$. However, $\,q\,$ is not parallel to $\,r_0 q\,$ if $\,r_0\in\mathbbm C.\,$ 
We notice that the orthogonality and parallelism concepts, respectively defined in (\ref{hn10}) and (\ref{hn12}), are identically valid for complex numbers. 
Let us consider further possibilities to the quaternionic basis (\ref{hn005}) and (\ref{hn007}). {\bf Remembering the definition of the general imaginary unit $\omega$ from (\ref{hn04}), one obtains}
\begin{equation}
 \big\langle q,\,\omega q\big\rangle=0,
\end{equation}
and  $q$ and $\omega q$ are consequently  orthogonal. In order to obtain the further two orthogonal components of the space, let us define:
\begin{definition}[Regular quaternionic curves] The parametrized curve $q:I\to\mathbbm H$ is regular in $\,I\subset\mathbbm R\,$ if $\;|q'(t)|\neq 0\quad\forall t\in I.$
\end{definition}
Thus, we can prove that
\begin{proposition}[ Imaginary quaternionic function ]\label{rnP00} Let $x=x(t)$ be a pure imaginary regular quaternionic curve. If 
$\, \omega=\nicefrac{x}{|x|},\,$ and $\,\omega'=\nicefrac{d\omega}{dt},\,$ then
\[
\begin{aligned}
\mbox{\bf i.} &\qquad\overline{\omega}'=-\,\omega'\\
\mbox{\bf ii.} &\qquad\omega\omega'=-\,\omega'\omega 
\end{aligned}
\]
\end{proposition}
{\bf Proof:} The item (i) comes immediately from the definition, and the item (ii) comes from the derivative of $\omega^2=-1$.

$\hfill\bm\square$

Consequently, {\bf assuming the derivative $\omega'$ of the imaginary unitary quaternion $\omega$ as a regular quaternionic function,} the four orthogonal directions defined from $q$ are as follows,
\begin{equation}
q\mapsto\big\{q,\,\omega q ,\,\omega'q,\,\omega\omega'q\big\},
\end{equation}
and the four-components of the orthogonal basis set comprises
\begin{equation}
\big\{1,\,i,\,j,\,k\big\}\equiv\left\{1,\,\omega ,\,\frac{\omega'}{|\omega'|},\,\omega\frac{\omega'}{|\omega'|}\right\}.
\end{equation}

In summary, we defined the quaternions as a four dimensional real vector space, and obtained several alternative basis for that. In the following section,
we will consider this subject using  a different notation.
\subsection{Polar notation\label{pn}}
Observing that the real coefficient of the imaginary component of (\ref{hn04}) is always positive, we choose
\begin{equation}
 x_0=\rho\cos\theta,\qquad |\,\vec x\,|=\rho\sin\theta
 ,\qquad\rho=|q|,\qquad\mbox{and}\qquad\theta\in[0,\,\pi].
\end{equation}
The polar notation of the Cartesian quaternion consequently comprises
\begin{equation}\label{hn05}
q\,=\,\rho\big(cos\theta+I\sin\theta\big),\qquad I=\cos\phi\, i+\sin\phi\, e^{i\xi}j,\qquad \phi\in[0,\,\pi],\quad\mbox{and}\qquad \xi\in[0,\,2\pi],
\end{equation}
Of course $\,I^2=-1,\,$ the inner product (\ref{hn005}) is valid, and $\,q\,$ and $\,Iq\,$ are orthogonal. In accordance to Proposition \ref{rnP00},  the two remaining orthogonal directions follow as derivatives of $\,I:$
\begin{equation}\label{hn07}
J=\frac{\partial{I}}{\partial \phi},\qquad K=\frac{1}{\sin\phi}\frac{\partial{I}}{\partial \xi},
\end{equation}
we obtain
\begin{equation}\label{hn08}
J=-\sin\phi i+\cos\phi e^{i\xi}j,\qquad\mbox{and}\qquad K=e^{i\xi}ij.
\end{equation}
The novel imaginary units satisfy the multiplication rule (\ref{hn003}) of the natural quaternionic imaginary units, and thus,
\begin{equation}\label{hn09}
\big\{1,\,i,\,j,\,k\big\}\,\equiv\,\big\{1,\,I,\,J,\,K\big\}.
\end{equation}
Consequently, we have the four orthogonal directions defined from $q$ in the polar notation as
\begin{equation}
q\mapsto\big\{q,\,I q ,\,J q,\,K q\big\}.
\end{equation}
Therefore, and in the same token as complexes, a perfect equivalence between the cartesian and the polar notation is observed in quaternions.

As a final remark, the range of $\theta\in[0,\,\pi]\,$ in (\ref{hn05}), and not $\,[0,\,2\pi]\,$ as in $\mathbbm C$,  indicates that a negative signal have to be absorbed by the quaternionic structure in order to keep the polar angles within the correct range.
Let us consider the simplest case:
\begin{equation}\label{hn14}
-q=(-1)\rho(\cos\theta+I\sin\theta).
\end{equation}
If $q\in\mathbbm C$, this multiplication means the rotation $\theta\to\theta+\pi$, what is prohibited if $q\in\mathbbm H$. Absorbing the negative signal in the imaginary unit, the correct rotation to quaternions will be 
\begin{equation}
 -q=\rho\Big[\cos(\pi-\theta)-I\sin(\pi-\theta)\Big],\qquad\mbox{and}\qquad
 -I=\cos(\pi-\phi)+\sin(\pi-\phi)e^{i(\xi+\pi)}j.
\end{equation}
Assuming $|q|=1$ for the sake of simplicity, in general we have
\begin{equation}\label{hn17}
q(\theta_1+\theta_2)\,=\,\left\{
\begin{array}{ll}
\cos\theta_0+\omega\sin\theta_0 & n\quad\mbox{even}\qquad \theta_0\in\left[0,\,\pi\right]\\ \\
\cos\left(\pi-\theta_0\right)-I\sin\left(\pi-\theta_0\right)& n\quad\mbox{odd,}\qquad \theta_0\in\left[0,\,\frac{\pi}{2}\right]\\ \\
\cos\left(\frac{\pi}{2}-\theta_0\right)-I\sin\left(\frac{\pi}{2}-\theta_0\right)
& n\quad\mbox{odd,}\qquad \theta_0\in\left[\frac{\pi}{2},\,\pi\right]\\ \\
\theta_1+\theta_2\,=\,n\pi+\theta_0.
\end{array}
\right.
\end{equation}
As a final aspect, the next subsection contains the symplectic notation for quaternions.
\subsection{Symplectic notation}
Quaternionic numbers can also be written as
\begin{equation}\label{sn01}
q=z_0+z_1 j,\qquad\mbox{where}\qquad z_0=x_0+x_1 i,\qquad\mbox{and}\qquad z_1=x_2 +x_3 i.
\end{equation}
The complex components of $q$ can be obtained from
\begin{equation}\label{sn02}
 z_0=\frac{1}{2}\big(q-iqi\big)\qquad\mbox{and}\qquad z_1=\frac{1}{2}\big(\,\overline q+i\,\overline q\,i\big)j.
\end{equation}
The symplectic notation is not unique, and we can replace (\ref{sn01}) with 
\begin{equation}\label{sn03}
q=z_0+\overline\zeta\, k \qquad\mbox{where}\qquad \zeta=x_3 +x_2 i,
\end{equation}
and four other possibilities replacing $i$ with $j$ and $k$. A symplectic polar notation reads 
\begin{equation}\label{sn04}
q\,=\,\rho\Big(\cos\vartheta e^{i\phi}+\sin\vartheta e^{i\psi}j\Big)\qquad\mbox{where}\qquad \vartheta\in\left[0,\frac{\pi}{2}\right]\qquad\mbox{and}
\qquad \phi,\,\psi\in[0,\,2\pi].
\end{equation}
In symplectic notation, the vector space is complex. Using (\ref{sn02}),  the inner product
\begin{equation}\label{sn05}
 \big\langle p,\,q\big\rangle=\frac{1}{2}\big(p\overline q-ip\overline qi\big),
\end{equation}
 satisfies 
\begin{equation}
\overline{\big\langle p,\,q\big\rangle}=\big\langle q,\,p\big\rangle,\qquad \big\langle q,\,q\big\rangle=|q|^2,\qquad
\big\langle \alpha_0 p,\,q\big\rangle=\alpha_0 \big\langle p,\,q\big\rangle,\qquad \big\langle p,\,\alpha_0 q\big\rangle=\big\langle p,\,q\big\rangle\overline\alpha_0,
\end{equation}
where $p,\,q\in\mathbbm H$ and $\alpha_0\in \mathbbm C$. The order within the product between the complex and the quaternion is inflexible, and $p\,\alpha_0$ is not allowed in this case. As in the cartesian case, the addition of the polar angle satisfies
\begin{equation}\label{hn19}
q(\vartheta_1+\vartheta_2)\,=\,
\left\{
\begin{array}{ll}
\cos\vartheta_0 e^{i\phi}+\sin\vartheta_0 e^{i\psi}j & \quad n=0\;\mbox{mod}\;4\\ \\
\cos\left(\frac{\pi}{2}-\vartheta_0\right)e^{i(\psi-\pi)}+\sin\left(\frac{\pi}{2}-\vartheta_0\right)e^{i\phi}j & \quad n=1\;\mbox{mod}\;4\\ \\
\cos\vartheta_0 e^{i(\psi-\pi)}+\sin\vartheta_0 e^{i(\phi-\pi)}j &\quad n=2\;\mbox{mod}\;4\\ \\
\cos\vartheta\left(\frac{\pi}{2}-\vartheta_0\right)e^{i\psi}+\sin\left(\frac{\pi}{2}-\vartheta_0\right)e^{i(\phi-\pi)}j&\quad n=3\;\mbox{mod}\;4,
\end{array}	
\right.
\end{equation}
where
\[
\vartheta_1+\vartheta_2=\vartheta_0 + n\frac{\pi}{2},\qquad\vartheta_0\in\left[0,\,\frac{\pi}{2}\right]
\qquad\mbox{and}\qquad n\in\mathbbm N.
\]
We finally have a comprehensive picture of quaternions to the purposes of differential geometry, introduced in the next section.

\subsection{Quaternionic gradient}
We adopt  the quaternionic gradient operator
\begin{equation}\label{c02}
\nabla =\sum_\mu \overline e_\mu\partial_\mu,
\end{equation}
that is already known from harmonic analysis \cite{Gilbert:1991ahi}, and that is different from the gradient vector $\bm\nabla$. Moreover, using the arbitrary quaternionic function
\begin{equation}\label{c03}
f=f^{(0)}+f^{(1)}i+f^{(2)}j+f^{(3)}k
\end{equation}
where  $\,f^{(\mu)}\,$ are real functions, let us entertain the equation
\begin{equation}\label{c04}
\nabla f=0\qquad\Rightarrow \qquad\sum_{\mu=0}^4\big(\nabla f\big)_\mu e_\mu=0
\end{equation}
which implies that each quaternionic component must be independently zero, so that
\begin{equation}\label{c05}
\left\{
\begin{array}{ll}
e_0: &f^{(0)}_0+f^{(1)}_1+f^{(2)}_2+f^{(3)}_3=0\\ \\
e_1: &f^{(1)}_0-f^{(0)}_1-f^{(3)}_2+f^{(2)}_3=0\\ \\
e_2: &f^{(2)}_0+f^{(3)}_1-f^{(0)}_2-f^{(1)}_3=0\\ \\
e_3: &f^{(3)}_0-f^{(2)}_1+f^{(1)}_2-f^{(0)}_3=0
\end{array}
\right.,
\qquad\mbox{where}\qquad
f^{(\mu)}_\nu=\partial_\nu f^{(\mu)}.
\end{equation}
 The conditions (\ref{c05}) are analogous to the Cauchy-Riemann conditions of complex analysis, although insufficent to define functions of a quaternionic variable. Taking the second derivative, so that
\begin{equation}
\overline\nabla\nabla f=0.
\end{equation}
we obtain
\begin{equation}
\nabla^2f^{(\mu)}=0,\qquad\mbox{where}\qquad\nabla^2=\bm{\nabla\cdot\nabla}
\end{equation}
is the usual Laplacian operator. This result is analogous to the complex one, and the real components of the quaternionic function $f$ that satisfies (\ref{c04}) are harmonic funtions in four dimensions, in the same fashion as harmonic functions of two dimensions solve the complex case. Therefore, the quaternionic case  generalizes  the complex one, as desired. Concluding this section, contemplate the quaternionic gradient in polar notation.
\subsection{Quaternionic gradient in polar notation}
Using the polar notation described in Section \ref{pn}, using the quaternionic basis (\ref{hn09}) we obtain   
\begin{eqnarray}\label{pc01}
\nonumber &i&\!\!\!\!=\cos\phi I-\sin\phi J\\
          &j&\!\!\!\!=\cos\xi\big(\sin\phi I+\cos\phi J\big)-\sin\xi K\\
\nonumber &k&\!\!\!\!=\sin\xi\big(\sin\phi I+\cos\phi J\big)+\cos\xi K,
\end{eqnarray} 
and  from polar expression of the $x_\mu$ cartesian coordinates, 
\begin{equation}\label{pc02}
\rho^2=\sum_{\mu=0}^3 x_\mu^2,\qquad \tan\theta=\frac{\sqrt{x_1^2+x^2_2+x^2_3}}{x_0},\qquad\tan\phi=\frac{\sqrt{x_2^2+x^2_3}}{x_2},\qquad
\tan\xi=\frac{x_3}{x_2},
\end{equation}
enable us to obtain the quaternionic gradient operator (\ref{c02}) in polar notation
\begin{equation}\label{pc03}
\nabla=\big(\cos\theta-I\sin\theta\big)\left(\partial_\rho-\frac{I}{\rho}\partial_\theta\right)-\frac{J}{\rho\sin\theta}\partial_\phi-\frac{K}{\rho\sin\theta\sin\phi}\partial_\xi.
\end{equation}
Using the quaternionic function
\begin{equation}\label{pc04}
g\big(\rho,\,\theta,\,\phi,\,\xi\big)=g^{(0)}+g^{(1)}I+g^{(2)}J+g^{(3)}K,
\end{equation}
in
\begin{equation}\label{pc05}
\nabla g=0
\end{equation}
 we obtain
\begin{equation}\label{pc06}
\left\{
\begin{array}{r}
\rho g^{(0)}_\rho+g^{(1)}_\theta+\cot\theta\left(2g^{(1)}+\cot\phi\, g^{(2)}+g^{(2)}_\phi+\frac{g^{(3)}_\xi}{\sin\phi}\right)
-\frac{g^{(2)}_\xi}{\sin\phi}+g^{(3)}_\phi+\cot\phi\,g^{(3)}=0\\ \\
\rho g^{(1)}_\rho-g^{(0)}_\theta+2g^{(1)}+\cot\phi\, g^{(2)}+g^{(2)}_\phi+\frac{g^{(3)}_\xi}{\sin\phi}+\cot\theta\left(
\frac{g^{(2)}_\xi}{\sin\phi}-g^{(3)}_\phi-\cot\phi\,g^{(3)}\right)=0\\ \\
\rho g^{(2)}_\rho+g^{(3)}_\theta-\cot\theta\left(g^{(0)}_\phi+\frac{g^{(1)}_\xi}{\sin\phi}-g^{(3)}\right)
+\frac{g^{(0)}_\xi}{\sin\phi}-g^{(1)}_\phi+g^{(2)}=0\\ \\
\rho g^{(3)}_\rho-g^{(2)}_\theta-g^{(0)}_\phi-\frac{g^{(1)}_\xi}{\sin\phi}+g^{(3)}
-\cot\theta\left(\frac{g^{(0)}_\xi}{\sin\phi}-g^{(1)}_\phi+g^{(2)}\right)=0,
\end{array}
\right.\qquad
\end{equation}
where we also employed the derivatives of the imaginary units
\begin{eqnarray}
&I_\phi=J,\qquad\qquad\;\, &J_\phi=-I,\qquad \qquad K_\phi=0,\\
&I_\xi=\sin\phi K,\qquad &J_\xi=\cos\phi K,\qquad \,K_\xi=\big(\cos\phi I-\sin\phi J\big)K.
\end{eqnarray}
Therefore, from
\begin{equation}
 \overline \nabla \nabla g=0,
\end{equation}
after several manipulations, we obtain
\begin{equation}
\nabla^2 g^{(\mu)}=0
\end{equation}
where the Laplacian operator reads
\begin{equation}
\nabla^2u=\frac{1}{\rho^3}\left(\rho^3 u_\rho\right)_\rho+\frac{1}{\rho^2\sin^2\theta}\left(\sin^2\theta\, u_\theta\right)_\theta+
\frac{1}{\rho^2\sin^2\theta\sin\phi}\left(\sin\phi u_\phi\right)_\phi+\frac{1}{\rho^2\sin^2\theta\sin^2\phi}u_{\xi\xi}.	
\end{equation}
Now, the description of quaternions is complete, and we can finally proceed to the description of the quaternionic differential geometry.

\section{QUATERNIONIC CONSTRAINTS\label{QC}} 
In a former article \cite{Giardino:2021onv}, we studied  quaternionic curves defining quaternionic Frenet-Serret equations. Taking benefit of this experience, and of Proposition \ref{rnP00}, in this section we generalize the formalism to higher dimensions. Let us start defining our simplest objects.

\begin{definition}[Quaternionic constraint]\label{qc01} A constraint in $\mathbbm H$ is the smooth map
$q(u,\,v,\,w): U\to\mathbbm H\,$ for $\,U\subset\mathbbm R^3\,$  such that
\[
q=x_0+x_1i+x_2 j+x_3k,
\]
 $\,x_\mu=x_\mu(u,\,v,\,w)\,$ are $\mathcal C^\infty$ real functions.
\end{definition}
In principle, the real functions $x_\mu$ may have an arbitrary number of variables. For the sake of simplicity, we choose only three, and left the more complicated cases to be entertained in future research. Therefore, we define:
\begin{definition}[tangent map]\label{qc02} Given a quaternionic constraint, the tangent map, and the unitary tangent map related to the parameter $a$ are such as 
\[
q_a=\frac{\partial q}{\partial a},\qquad \mbox{and}\qquad t^{(a)}=\frac{q_a}{|q_a|}\qquad\mbox{where}\qquad a=\big\{u,\,v,\,w\big\}.
\]
\end{definition}
It is implicit that the derivative operates over the real components of the quaternionic function. A precise definition to quaternionic constraints would be
\begin{definition}[Regular quaternionic subspace] \label{qc03}
A quaternionic subset $\,S\subset\mathbbm H\,$ is regular if $\,\forall\,q_0\in S\,$ there are an open neighborhood $\,V\subset\mathbbm H,\,$ a subset $\,U\in\mathbbm R^n\,$ and a bijection $q:U\to V\cap S$ such that
\begin{enumerate}
\item $n\in\big\{1,\,2,\,3\big\}$.
\item the $x_\mu\,$ components of $q$ belong to the $\,\mathcal C^\infty\,$ class.
\item $x_\mu$ are  homeomorphisms.
\item $\,\forall\,P_0\in U\,$ the first order partial derivatives $q_a(P_0)$ are non zero and linearly independent.
\end{enumerate}
\end{definition}
A regular quaternionic contraint $q$ defines  a quaternionic subspace, emphasizing that $q$ is not a quaternionic variable function. This quaternionic subspace can be a quaternionic curve, a quaternionic surface, or a quaternionic three-dimensional hyper-surface. Let us impose unitary tangent vectors, where $\,\big|t^{(a)}\big|=1.\,$  In this case, it holds the orthogonality condition
\begin{equation}\label{qc04}
\left\langle t^{(a)},\,t^{(a)}_b\right\rangle = 0,
\end{equation}
where $t^{(a)}_b\,$ is the derivative of $t^{(a)}$ with respect to $b$. Since $e_m t^{(a)}$ is orthogonal to $t^{(a)}$,   the second derivative  belongs to a quaternionic subspace generated from the product of the tangent component and the imaginary unit. Consequently,
\begin{equation}\label{qc05}
\left(\frac{q_{a}}{|q_a|}\right)_a=\kappa^{(a)}q_a\qquad\mbox{where}\qquad \kappa^{(a)}=\sum_{\ell=1}^3\kappa^{(a)}_\ell e_\ell,
\end{equation}
and $\kappa^{(a)}$ is called the quaternionic curvature function. Additionally, 
\begin{equation}\label{qc06}
\left(\frac{q_{a}}{|q_a|}\right)_b=\tau^{(ab)}q_a\qquad\mbox{where}\qquad \tau^{(ab)}=\sum_{\ell=1}^3\tau^{(ab)}_\ell e_\ell,
\end{equation}
and $\tau^{(ab)}$ is the quaternionic torsion function. The evident analogy to the Frenet-Serret equations of real curves enables us to define the unitary quaternionic normal $n^{(a)}$, and the unitary quaternionic binormal $b^{(ab)}$, so that 
\begin{equation}\label{qc07}
n^{(a)}=\frac{1}{|\kappa^{(a)}|}\kappa^{(a)}t^{(a)},\qquad\qquad\mbox{and}\qquad\qquad b^{(ab)}=\frac{1}{|\tau^{(ab)}|}\tau^{(ab)}t^{(a)}.
\end{equation}
The real components of the curvature are
\begin{equation}\label{qc08}
\kappa^{(a)}_\ell=\frac{1}{|q_a|^2}\big\langle q_{aa},\,e_\ell t^{(a)}\big\rangle
\end{equation}
and the real components of the torsion are
\begin{equation}\label{qc09}
\tau^{(ab)}_\ell=\frac{1}{|q_a|^2}\big\langle q_{ab},\,e_\ell t^{(a)}\big\rangle,
\end{equation}
and we immediately observe that, in principle, $\,\tau^{(ab)}\neq\tau^{(ba)}.\,$  We remark that (\ref{qc05}-\ref{qc06}) have in fact a second possibility, because we could have written the curvature and the torsion at the right hand side of the quaternionic funtion. Thus, our definitions are in fact a convenience that can be modified. For example, if we rotate our function using a constant unitary quaternion $u$, we have
\begin{equation}
\left(\frac{uq_{a}}{|q_a|}\right)_a=\widetilde\kappa^{(a)}uq_a,\qquad\mbox{where}\qquad \widetilde\kappa^{(a)}=u\kappa^{(a)}\overline u.
\end{equation}
However, using the transformation $qu$, it would be more convenient to use $q_a\,\kappa^{(a)}$ on the right hand side of (\ref{qc05}).
Let us prove a final and interesting result of this section:
\begin{lemma}\label{l01}
If $q$ is a regular quaternionic curve, we have
\end{lemma}
\begin{eqnarray}    
\mbox{\bf i.}\nonumber && q_{aa}\,=\,\big|q_a\big|_a\,t^{(a)}\,+\,\big|q_a\big|^2\big|\kappa^{(a)}\big|\,n^{(a)}\\
\mbox{\bf ii.}\nonumber && q_{ab}\,=\,\big|q_a\big|_b\, t^{(a)}\,+\,\big|q_a\big|^2\big|\tau^{(ab)}\big|\,b^{(ab)}\,=\,
\big|q_b\big|_a t^{(b)}\,+\,\big|q_b\big|^2\big|\tau^{(ba)}\big|b^{(ba)}.
\end{eqnarray}
{\bf Proof:} The above relations are immediately obtained differentiating $q_a=|q_a|t^{(a)}$ respectively to $a$ and $b$ and using (\ref{qc05}-\ref{qc07}).

$\hfill\bm\square$

We can interpret Lemma \ref{l01}-i in dynamical terms, where $q_{aa}$ is the acceleration of a point of the quaternionic curve of velocity $q_a$. The normal component of the acceleration comprises the centripetal acceleration, whose radius of the curvature is
\begin{equation}\label{qc10}
R^{(a)}=\frac{1}{\big|\kappa^{(a)}\big|}.
\end{equation}
And the analogy to Newtonian dynamics is exact. By analogy, we define the torsion radius
\begin{equation}\label{qc13}
R^{(ab)}=\frac{1}{\big|\tau^{(ab)}\big|}.
\end{equation}
and the interpretation is analogous.

\section{QUATERNIONIC DIFFERENTIAL FORMS}

The quaternionic directional derivative is defined to be
\begin{definition}[quaternionic directional derivative] Let $f:\mathbbm R^4\to\mathbbm H$ be a quaternionic regular function, $\,\nabla\,$ the quaternionic gradient operator, and $q$ an arbitrary quaternion number. The inner product
\[
D_qf(P_0)=\big\langle q, \nabla f\big\rangle\Big|_{P_0}
\]
defines the quaternionic directional derivative $D_q f:\mathbbm R^4\to \mathbbm R$, along the direction of $q$, and at the point $P_0\in\mathbbm R^4$.
\end{definition}

The above definition can be interpreted as a total real derivative of $f$ in terms  a parameter $t$. Recalling the correspondence established by the quaternionic gradient between the coordinates $x_\mu$ and $e_\mu$, if $x_\mu=p_{0\mu}+tq_\mu$, we can use the quaternion $q=\sum e_\mu q_\mu$ in the directional derivative in the same fashion as the real directional derivative in terms of a parameter $t$. Therefore, familiar properties of the real directional derivative can be proven, namely
\begin{proposition} Let $f$ and $g$ be quaternionic functions on $\mathbbm R^4$, and the constants $p,\,q\in\mathbbm H$, and $\alpha,\,\beta\in\mathbbm R$. Then
\[
\left\{
\begin{array}{ll}
\mbox{\bf i.} & \alpha D_p f+\beta D_q f=\Big(\alpha D_p +\beta D_q \Big)f\\ \\
\mbox{\bf ii.} & D_p\big(\alpha f+\beta g)=\alpha D_pf+\beta D_pg\\ \\
\mbox{\bf iii.} & D_p\big(fg\big)=\big\langle(\nabla f)g,\,p\big\rangle+\big\langle f(\nabla g),\,p\big\rangle \\
\end{array}
\right.
\]
\end{proposition}
The proofs of ({\bf i}) and ({\bf ii}) are immediate, while ({\bf iii}) demands to verify the Leibniz rule to the quaternionic gradient,
\[
 \nabla(fg)=(\nabla f)g+f\nabla g,
\]
which is also straightforward.

$\hfill\square$

By way of example,  the quaternionic directional derivative of the cartesian coordinates multiplied by the quaternionic basis elements gives
\begin{equation}\label{qdd01}
D_px_{\mu\nu}=
\left\{
\begin{array}{l}
D_p (x_0e_0)=p_0\\
D_p (x_me_0)=-p_m\\
D_p (x_0e_m)=p_m\\
D_p (x_me_n)=\delta_{mn} p_0+\epsilon_{mn\ell}p_\ell
\end{array}
\right..
\end{equation}
In  matrix form, we can organize the real components
\begin{equation}\label{qdd02}
D_px_{\mu\nu}=
\left[\begin{array}{rrrr}
 p_0 & p_1 & p_2 & p_3\\
 -p_1 & p_0 & p_3 & -p_2\\ 
 -p_2 & -p_3 & p_0 & p_1\\
 -p_3 & p_2 & -p_1 & p_0
  \end{array}
\right].
\end{equation}
Here we observe an exact concordance between (\ref{qdd01}-\ref{qdd02}) and (\ref{c04}-\ref{c05}), and therefore, we define the quaternionic differential
\begin{definition}[Quaternionic differential $1-$form]\label{def_form} The differential $1-$form associated to the map $f:\mathbbm R^4\to\mathbbm H$ is
\begin{align*}
 df&=\sum_{\mu,\,\nu=0}^3f_{\mu\nu}Dx_{\mu\nu}\\
 &=\sum_{\mu=0}^3\big(\nabla f\big)_\mu dx_\mu.
\end{align*}
\end{definition}
The agreement to the usual real differential forms is exact, and the usual properties of one forms, namely their anti-commutativity property, is immediately preserved within this formulation. Moreover, the wedge product of $1-$forms is also immediately obtained.

\section{CONNECTION}
Let us consider a unitary quaternion $\,u\,$ that can be deployed  as a basis to quaternions following (\ref{hn007}). Therefore, in the same fashion as (\ref{qc05}), we can write
\begin{equation}
 u_a=\omega^{(a)} u,\qquad\Rightarrow\qquad \omega^{(a)}=u_a\overline u,
\end{equation}
where $\omega^{(a)}$ is of course a pure imaginary quaternion. A similar situation compared to the usual connection $1-$form of real differential geometry,  observing that (\ref{qc05}) can be interpreted as a quaternionic covariant derivative. In matrix terms, we observe that $\omega^{(a)}$ is accordingly anti-symmetric, and their real components are 
\begin{eqnarray}\label{cd03}
\omega^{(a)}_\ell=\big\langle\, u_a,\,e_\ell u\,\big\rangle.
\end{eqnarray}
Moreover, the perfect correspondence between the quaternionic $\omega^{(a)}$  and the real connection in matrix formalism enables us to 
interpret $\omega^{(a)}$ as the quaternionic connection. Using this fact, the fundamental structure of real differential geometry comes immediately, providing
the translation of the quaternionic formalism into the language of the real differential geometry. Hence, we define
\begin{definition}[Quaternionic dual $1-$forms] If $e_\mu$ is a basis for $\mathbbm H$ , then $,\forall\,p\in\mathbbm H$ the dual $1-$form is defined by
\[
 \phi_\mu(p)=\langle p,\,e_\mu u\rangle.
\]
\end{definition}
In the case of the natural basis (\ref{hn005}),
\begin{equation}
 dx_\mu(p)=\langle p,\,e_\mu \rangle.
\end{equation}
 We also observe that the components of the general quaternionic basis (\ref{hn007}) written as
 \begin{equation}
  q_\mu=e_\mu q
 \end{equation}
generate
\begin{equation}
 \phi_\mu=\sum_\nu a_{\mu\nu}dx_\mu
\end{equation}
where the real coefficient $a_{\mu\nu}$ are the real components of the quaternionic basis element $q_\mu$. The sole difference to the real case is the constrained standard of $a_{\mu\nu}$, that obeys (\ref{qdd02}). Therefore, in complete analogy to real differential geometry, we have the following
lemma that is a direct consequence of the linearity of $1-$forms.
\begin{lemma}
Let $\,\phi_\mu\,$ be the $1-$form basis dual to $\,g_\mu.\,$ Thus an arbitrary $1-$form $\,\psi\,$ has the unique expression
\[
\psi=\sum_\mu \psi(g_\mu)\phi_\mu.
\]
\end{lemma}

Therefore, also in analogy to the established knowledge of differential geometry, we conclude that the structural equations hold, namely
\begin{equation}
d\phi=\omega^{(a)}\wedge \phi,\qquad\qquad
d\omega=\omega^{(a)}\wedge\omega^{(a)}.
\end{equation}
where $\phi$ is the $1-form$ dual to the natural quaternionic basis, and where the $1-$forms are in matrix representation. We must also stress that the exterior derivative is such as Definition \ref{def_form}. As a final remark of this section, we obtained a translation of the quaternionic structure into the real structure by using the matrix representation, and this simple result is the structural fact that allows the establishment of $\mathbbm H$DG.

\section{FINAL CONSIDERATIONS}
In this article we provided the foundations of a quaternionic differential geometry ($\mathbbm H$DG), where  quaternionic counterparts of the structural elements of the differential geometry of real curves and surfaces were obtained. There are several possible interpretations to this new geometry. First of all, it enables to build pure quaternionic figures within a quaternionic space, and hence we have defined a geometry of quaternionic objects in analogy to the real geometric figures. 

On the other hand, we may have a geometry that is not analogous to the geometry of real objects. This can be observed mainly in the components of the quaternionic differential $1-$forms, which are not identically related to the real $1-$forms, and therefore the calculus that will arise from this novel geometry must be fundamentally different from that obtained from the real one. This is possibly the most interesting fact risen up from this novel geometry, and whose consequences have to be investigated in future investigation. There are several possibilities for that, such as the calculus, but also concerning the Riemannian version of $\mathbbm H$DG. We hope that in a close future we can answer these important questions, and open a broad field of investigation of differential geometry and their physical applications.

\paragraph{Funding} This work is supported by the Funda\c c\~ao de Amparo \`a Pesquisa do Rio Grande do Sul, FAPERGS, within the Edital 14/2022.

\paragraph{Conflict of interest statement}The author declares that he has no known competing financial interests or personal relationships that
could have appeared to influence the work reported in this paper.

\paragraph{Data availability statement}The author declares that data sharing is not applicable to this article as no datasets were generated or analysed during the current study.

%
%
%
%
\begin{footnotesize}

\begin{thebibliography}{10}

\bibitem{Giardino:2021onv}
Giardino, S.:
\newblock {\it A primer on the differential geometry of quaternionic curves}.
\newblock  Math. Methods Appl. Sci.  44 (18):14428--14436 (2021).  https://doi.org/10.1002/mma.7709

\bibitem{Bharathi:1987kbn}
Bharathi, K., Nagaraj, M.:
\newblock {\it Quaternion valued function of a real Variable Serret-Frenet formulae}.
\newblock  Indian J. Pure Appl. Math.  18:507--511 (1987).

\bibitem{Gunes:1994mmt}
Sivridag, A. I., Gunes, R., Keles, S.:
\newblock {\it The Serret-Frenet formulae for dual quaternion-valued functions of a single real variable}.
\newblock Mech. Mach. Theor.  29(5):749--754 (1994). https://doi.org/10.1016/0094-114X(94)90116-3

\bibitem{Girard:2015lnc}
Girard, P. R., Clarysse, P., Pujol, R., Wang, L., Delachartre, P.:
\newblock {\it Differential Geometry Revisited by Biquaternion Clifford Algebra}.
\newblock  In: J. D. Boissonnat et al. (eds) Curves and Surfaces 2014.
  Lecture Notes in Computer Science  9213 Springer, Cham(2):47--64  (2015).

\bibitem{Aksoyak:2019}
Aksoyak, F. K.:
\newblock {\it A new type of quaternionic frame in $\mathbbm R^4$}.
\newblock Int. J. Geom. Meth. Mod. Phys.  16(6):1959984 (2019). https://doi.org/10.1142/S0219887819500841

\bibitem{Coken:2004otq}
Coken, A. C., Tuna, A.:
\newblock {\it On the quaternionic inclined curves in the semi-Euclidean space $E^4_2$''}.
\newblock Appl. Math. Comput. A155(2):373--389 (2004). https://doi.org/10.1016/S0096-3003(03)00783-5

\bibitem{Hacisalihoglu:2011gok}
G\"ok, I., Okuyucu, O. Z., Kahraman, F., H. H. Hacisalihoglu, H. H.:
\newblock {\it On the quaternionic $B_2$-slant helices in the Euclidean space $E^4$}.
\newblock Adv. Appl. Clifford Algebras  21:707--719,(2011). https://doi.org/10.1007/s00006-011-0284-6

\bibitem{Gungor:2011scq}
Gungor, m. A., Tosun, M.:
\newblock {\it Some characterizations of quaternionic rectifying curves}.
\newblock Differ. Geom. Dyn. Syst. 13:89--100 (2011).

\bibitem{Ilarslan:2013bma}
Kecilioglu, O., Ilarslan, K.:
\newblock {\it Quaternionic Bertrand curves in Euclidian $4-$space}.
\newblock Bull. Math. Anal. Appl. 5(3):27--38 (2013).

\bibitem{EJPAM1979}
Bektas, O., Gurses, N. B., Yuce, S.:
\newblock {\it Osculating Spheres of a Semi Real Quaternionic Curve in $E_2^4$}.
\newblock Eur. J. Pure and Appl. Math. 7(1):86--96 (2014).

\bibitem{Bektas:2016qoc}
Bektas, O., Gurses, N., Yuce, A.:
\newblock {\it Quaternionic osculating curves in Euclidean and semi-Euclidean space}.
\newblock J. Dyn. Sys. Geom. Theor.  14(1):65--84 (2016). https://doi.org/10.1080/1726037X.2016.1177935

\bibitem{Kisi:2017aap}
Ozturk, G., Kisi, I., Buyukkutuk, S.:
\newblock {\it Constant ratio quaternionic curves in Euclidean spaces}.
\newblock Adv. Appl. Clifford Algebras, 27:1659--1673 (2017). https://doi.org/10.1007/s00006-016-0716-4

\bibitem{Coken:2017nqc}
Coken, A. C., Tuna Aksoy, A.:
\newblock {\it Null quaternionic Cartan helices in $\mathbbm R^3_v$}.
\newblock Acta. Phys. Pol. A132(3-II):896--899 (2017). https://10.12693/APhysPolA.132.896

\bibitem{Karadag:2019scq}
Karadag, M., Sivridag, A. I.:
\newblock {\it Some characterizations for a quaternion-valued and dual variable curve}.
\newblock  Symmetry 11(2):125 (2019). https://doi.org/10.3390/sym11020125

\bibitem{Kizilay:2020eqc}
Kizilay, A., Yildiz, O. G., Okuyucu, O. Z.:
\newblock {\it Evolution of quaternionic curve in the semi-Euclidean space $E_2^4$}.
\newblock Math. Meth. Appl. Sci. 44(9):7577-7587 (2021). https://doi.org/10.1002/mma.6374

\bibitem{Kahraman:2020dqn}
Kahraman, T.:
\newblock {\it Differential equations of null quaternionic curves}.
\newblock Int. J. Appl. Comput. Math. 6(63):6583--6592 (2020). https://doi.org/10.1007/s40819-020-00824-3

\bibitem{Gungor:2013qie}
Soyfidan, T., Gungor, M. A.:
\newblock {\it On the quaternionic involute-evolute curves''}.
\newblock Preprint arXiv:1311.0621[math.GT] (2013).

\bibitem{Onder:2020gqi}
Hanif, M., \"Onder, M.:
\newblock {\it Generalized quaternionic involute-evolute curves in the Euclidean four-space $E^4$}.
\newblock Math. Meth. Appl. Sci. 43(7):4769--4780 (2020). https://doi.org/10.1002/mma.6231

\bibitem{Altun:2017scy}
Senyurt, S., Cevahir, C., Altun, Y,:
\newblock {\it On spatial quaternionic involute curve: a new view}.
\newblock Adv. Appl. Clifford Algebras 18:1815--1824 (2017). https://doi.org/10.1007/s00006-016-0669-7

\bibitem{Senyurt:2017scy}
Senyurt, S., Cevahir, C., Altun, Y,:
\newblock {\it On the Smarandache curves of spatial quaternionic involute curve}.
\newblock Proc. Natl. Acad. Sci. India A Phys. Sci.:1815--1824 (2019). https://doi.org/10.1007/s40010-019-00640-5

\bibitem{Hou:2018gia}
Hanif, M., Hou, Z. H.:
\newblock {\it Generalized involute and evolute curve-couple in Euclidean space}.
\newblock Int. J. Open Problems Compt. Math. 11(2):28--39 (2018).

\bibitem{Yayli:2016sqc}
Aslan, S., Yayli, Y.:
\newblock {\it Split quaternions and canal surfaces in Minkowski $3-$space}.
\newblock Int. J. Geom. 5(2):51--61 (2016).

\bibitem{Yayli:2016csq}
Aslan, S., Yayli, Y.:
\newblock {\it Canal surfaces with quaternions}.
\newblock Adv. Appl. Clifford Algebras 26(2):31--38 (2016). https://doi.org/10.1007/s00006-015-0602-5

\bibitem{Yayli:2017qso}
Aslan, S., Yayli, Y.:
\newblock {\it Quaternionic shape operator}.
\newblock Adv. Appl. Clifford Algebras  27(2):2921--2931 (2017). https://doi.org/10.1007/s00006-017-0804-0

\bibitem{Gok:2017sni}
G\"ok, I.:
\newblock {\it Quaternionic approach of canal surfaces constructed by some new ideas}.
\newblock Adv. Appl. Clifford Algebras 27(2):1175--1190 (2017). https://doi.org/10.1007/s00006-016-0703-9

\bibitem{Kocakusakli:2017cks}
Kocakusakli, E., Tuncer, O., G\"ok, I., Yayli, Y.:
\newblock {\it A new representation of canal surfaces with split quaternions in Minkowski $3-$Space}.
\newblock Adv. Appl. Clifford Algebras 27:1387--1409 (2017). https://doi.org/10.1007/s00006-016-0723-5

\bibitem{Karakus:2019ame}
Karakus, S. O.:
\newblock {\it Quaternionic approach on constant angle surfaces in $S^2\times\mathbbm R^2$}.
\newblock Appl. Math. e-not. 19:497--506 (2019).

\bibitem{Yayli:2019zco}
Canakci, Z., Tuncer, O. O., G\"ok, I., Y. Yayli, Y.:
\newblock {\it The construction of circular surfaces with quaternions}.
\newblock Asian-Eur. J. Math. 12(1):1950091 (2019). https://doi.org/10.1142/S1793557119500918

\bibitem{Bekar:2021rsc}
Aslan, S., Bekar, M., Yayli, Y.:
\newblock {\it Ruled surfaces constructed by quaternions}.
\newblock J. Geom. Phys. 161:104048 (2021). https://doi.org/10.1016/j.geomphys.2020.104048

\bibitem{Tuncer:2022mgt}
Tuncer, O. O.:
\newblock {\it Generalized tubes in pseudo-Galilean $3-$space: Split semi-quaternionic representations and an application to magnetic flux tubes}.
\newblock Math. Meth. Appl. Sci. 45(3):1468--1487 (2022). https://doi.org/10.1002/mma.7866

\bibitem{Berndt:1991jra}
Berndt, J.:
\newblock {\it Real hypersurfaces in quaternionic space forms}.
\newblock Journal f\"ur die reine und angewandte Mathematik 419(2):9--26 (1991). https://doi.org/10.1515/crll.1991.419.9

\bibitem{Suh:1997dga}
Perez, J. D., Suh, Y. J.:
\newblock {\it  Real hypersurfaces of quaternionic projective space satisfying $\nabla_{U_i}R=0$}.
\newblock Diff. Geom. Appl. 7(3):211--217 (1997). https://doi.org/10.1016/S0926-2245(97)00003-X

\bibitem{Gori:2021ana}
Gentili, G., Gori, A., Sarfatti, G.:
\newblock {\it On compact affine curves and surfaces}.
\newblock J. Geom. Anal. 31:1073--1092 (2021). https://doi.org/10.1007/s12220-019-00311-2

\bibitem{Ward:1997qcn}
Ward, J. P.:
\newblock  Quaternions and Cayley Numbers.
\newblock  Springer Dordrecht (1997).

\bibitem{Garling:2011zz}
Garling, D. J. H.:
\newblock Clifford algebras: an introduction.
\newblock Cambridge Univ. Press (2011).

\bibitem{Rocha:2013qtt}
Vaz, J., da Rocha, R.:
\newblock An introduction to Clifford algebras and spinors.
\newblock Oxford University Press (2016).

\bibitem{Morais:2014rqc}
Morais, J. P., Georgiev, S., Spr\"ossig, W.:
\newblock Real quaternionic calculus handbook.
\newblock {Birkh\"auser} (2014).

\bibitem{Hurwitz:1898hqt}
Hurwitz, A.:
\newblock {\it Ueber die Composition der quadratischen Formen von belibig vielen Variablen}.
\newblock Nachr. Gesell. Wiss. G\"ottingen, Math-Phys. Kl. 309-316 (1898).

\bibitem{Harvey:1990sca}
Reese Harvey, F.:
\newblock Spinors and calibrations.
\newblock Academic Press (1990).

\bibitem{Gilbert:1991ahi}
Gilbert. J. E., M. A. M. Murray, M. A. M.:
\newblock Clifford algebras and Dirac operators in harmonic analysis.
\newblock Cambridge Univ. Press (1991).

\end{thebibliography}

\end{footnotesize}

\end{document}